\documentclass[twoside,leqno]{article}

\usepackage[letterpaper]{geometry}

\usepackage{siamproceedings}

\usepackage[T1]{fontenc}
\usepackage{amsfonts}
\usepackage{graphicx}
\usepackage{epstopdf}
\usepackage{enumitem}
\usepackage{algorithmic}
\ifpdf
  \DeclareGraphicsExtensions{.eps,.pdf,.png,.jpg}
\else
  \DeclareGraphicsExtensions{.eps}
\fi

\DeclareTextFontCommand{\emph}{\em\bf}

\newsiamremark{remark}{Remark}
\newsiamremark{hypothesis}{Hypothesis}
\crefname{hypothesis}{Hypothesis}{Hypotheses}
\newsiamthm{claim}{Claim}
\newsiamthm{problem}{Problem}
\newsiamthm{conjecture}{Conjecture}
\newsiamthm{question}{Question}
\def\cE{{\cal E}}
\def\cM{{\cal M}}
\def\cV{{\cal V}}
\def\cF{{\cal F}}
\def\cG{{\cal G}}
\def\cB{{\cal B}}
\def\cQ{{\cal Q}}
\def\cN{{\cal N}}
\def\se{{\mathsf e}}
\def\A{\mathbb A}
\def\B{\mathbb B}
\def\C{\mathbb C}
\def\R{\mathbb R}

\def\T{\mathbb T}
\def\D{\mathbb D}
\def\bS{\mathbb S}

\def\Z{\mathbb Z}
\def\cU{\mathcal U}
\def\cF{\mathcal F}
\def\cG{\mathcal G}

\def\leb{\mathrm{Leb}\, }
\def\Leb{\mathrm{Leb}\, }
\def \cP{{\mathscr P}}

\def\Sym{\mathrm{Symp}}
\def\Symp{\mathrm{Symp}}
\def\Diff{\mathrm{Diff}}
\usepackage{mathrsfs}
\def\sg{{\mathsf g}}
\def\sG{{\mathsf G}}

\begin{document}

\title{\Large Wild dynamics on manifolds}
    \author{  Pierre Berger\thanks{IMJ-PRG, CNRS, Sorbonne Université, Université Paris Cité, partially supported by the ERC project 818737 Emergence of wild differentiable dynamical systems.}
   }

\date{}

\maketitle

\begin{abstract} We survey a few  results  on differentiable, symplectic, or analytic wild dynamics.
\end{abstract}

\begin{center}
 \includegraphics[height=7.cm]{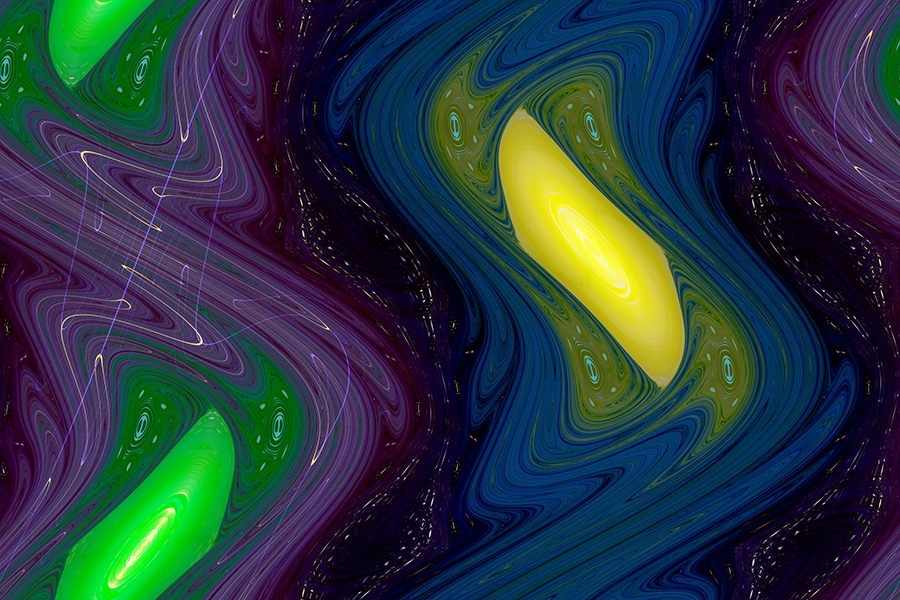}\\ 
 {A wild dynamics on the torus.}  
\end{center}

\section{Wild hyperbolicity}
The discovery of wild dynamics goes back to Poincaré~\cite{Po92}, when he realized the geometric complexity of  stable and unstable manifolds intersecting in a non-trivial way:
\begin{center} \em
`` ces intersections forment une sorte de treillis, de tissu, de réseau à mailles infiniment serrées ; chacune de ces courbes ne doit jamais se recouper elle-même, mais elle doit se replier elle-même d’une manière très complexe pour venir couper une infinité de fois toutes les mailles du réseau. ''
\end{center} 
In the late 60s, there were some hopes that the uniformly hyperbolic theory would  enable us to describe
most of differentiable systems~\cite[\S 1.6]{Sm67}, in particular the geometry of  their stable and unstable manifolds.
Soon it became clear that the mathematical reality was actually  much more complex.
\subsection{Definition and examples} 
The first counterexample showing that uniform hyperbolicity is not dense was found by Abraham and Smale~\cite{AS70} for diffeomorphisms in dimension 4 and then by Newhouse for  surface diffeomorphisms~\cite{Ne74}. The latter called his counterexample a \emph{wild hyperbolic set}:
\begin{definition} 
A hyperbolic compact set $K$ of a diffeomorphism $f$ is \emph{wild} if it displays local stable and unstable laminations $W^s_{loc} (K; f)$ and $W^u_{loc}(K; f)$ which are not transverse (at least one leaf of $W^s_{loc} (K; f)$ is tangent to one leaf of $W^u_{loc}(K; f))$, and this tangency is \emph{robust}:  for every
perturbation $\tilde f$ in an open  neighborhood $\cN$ of $f$, the 
 continuations $W^s_{loc} (K; \tilde f)$ and $W^u_{loc}(K; \tilde  f)$ are not transverse. The open set $\cN$, when taken maximal, is called the \emph{Newhouse domain}. 
\end{definition} 
Newhouse's initial article provides an example of  smooth surface diffeomorphism displaying a wild hyperbolic set, for perturbations in any $C^r$-topology, $\infty\ge r\ge 2$.
In~\cite{Ne79}, he proved that wild hyperbolic sets appear whenever we unfold a quadratic homoclinic tangency of a dissipative saddle point. As an application he obtained that wild hyperbolic sets appear in the quadratic real dissipative Hénon family. The counterpart of~\cite{Ne79} for surface symplectomorphisms has been shown by Duarte~\cite{Du08}.   
    Examples of wild hyperbolic sets have been found in many other dynamical  spaces, let us mention the Chirikov standard map family~\cite{Du94}, 
the space of complex polynomial automorphisms~\cite{Bu97},    the space of convex billiards~\cite{Ca22}, and the space of steady Euler flows~\cite{BFP23}.

The mechanism revealed by Abraham-Smale is called a \emph{robust heterodimensional cycle}. The construction used a hyperbolic set which   is now encompassed in the class of  \emph{blenders}~\cite{BDV05}:
 \begin{definition}  
 A \emph{blender}  $K$ for a diffeomorphism $f$  is a hyperbolic basic set, with stable and unstable dimensions $d_s$ and $d_u$,  such that there exists a local unstable lamination $W^u_{loc} (K;f)$ 
and    a submersion $\pi: N\to \R^{d_u+1}$  from a neighborhood $N$ of $K$ such that  $\pi( W^u_{loc} (K;f))$ contains a nonempty open subset $O$. Moreover,  for every sufficiently small $C^1$-perturbation $\tilde f$ of $f$,  a continuation $W^u_{loc} (K;\tilde f)$  has its image under $\pi$ which contains~$O$. 
\end{definition} 
\begin{remark} 
We notice that when $f$ is a local diffeomorphism of a surface $M$,  a hyperbolic basic set $K$ is a blender
if it has a local unstable set which robustly contains an open subset of $M$.
% if it displays a local unstable set $W^u_{loc} (K;f)$ which contains robustly an open subset of $ M$.  
\end{remark} 

There is a vast literature on the notion of wild hyperbolic set (and other topics discussed below), which we do not attempt to survey completely.

\subsection{Infinitude of the number of attractors} 
A long-standing problem is that of ergodicity of a typical dynamical system. This goes back to the Boltzmann ergodic hypothesis, which has been reformulated in modern terms by Birkhoff-Koopman~\cite{BK32}   as follows: a typical proper Hamiltonian system is ergodic on a.e. [component of] energy level. Here ergodicity means that   Liouville a.e. point has its  orbit which is equidistributed over the Liouville probability measure. This was disproved by Kolmogorov in 1954, with the celebrated  KAM theorem.

A great idea of Smale was to remove and simplify the structure left invariant by the system. Namely, he proposed to focus on differentiable systems -- which do not need to preserve the volume -- on low-dimensional and  compact manifolds.  He conjectured the open-density of uniform hyperbolicity in the space of such  dynamical systems. As a consequence of the   works of  Sinai and Bowen-Ruelle, the Smale conjecture would have implied the open-density of  \emph{almost finite ergodicity}:  the existence of a finite set of probability measures 
which model the statistical behavior of the orbit of Leb. a.e. point. 

From this perspective, the most disturbing point of Newhouse's initial article~\cite{Ne74} is that,  whenever a wild hyperbolic set $K$ of a surface diffeomorphism displays an area contracting periodic point, then  a topologically generic\footnote{
A \emph{(topologically) generic set} is a countable intersection of open-dense subsets.} perturbation of the dynamics displays infinitely many sinks\footnote{A \emph{sink} is an attracting periodic orbit.} which accumulate on the hyperbolic set.  Consequently, a  topologically generic perturbation is \emph{not} almost finitely ergodic.  Then Newhouse~\cite{Ne79}  showed that  the  coexistence of infinitely many attractors also occurs at generic parameters of any unfolding of a quadratic homoclinic tangency of an area contracting periodic point.  The same phenomenon was shown to occur in the space of complex Hénon maps by  Buzzard~\cite{Bu97} and in the $C^1$-topology by Bonatti and D\'iaz~\cite{BD99}. The counterpart of this result for surface conservative maps, was shown by Duarte~\cite{Du08}, where the sinks are replaced by elliptic islands (elliptic periodic points surrounded by arbitrarily small  KAM circles). 
Yet a topologically generic set can have empty interior   and, when it is in an Euclidean space $\R^n$,  it might be of Lebesgue measure zero.  In particular the following problem\footnote{
This problem has been formulated by many as a conjecture, in particular by  Palis~\cite{Pa00}.} is still open:
\begin{problem}  \label{Density finite}
Is there a dense subset of $\Diff(M)$ formed by dynamics with finitely many attractors, whose basin union covers $\leb$ a.e. $ M$? 
\end{problem}

\begin{figure}[h]
\begin{center} 
 \includegraphics[height=4.cm]{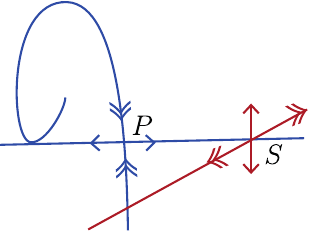} 
 \caption{The \emph{bicycle} configuration: a dissipative saddle point $P$ displaying both a homoclinic tangency and a 
 heterodimensional cycle with a   projectively hyperbolic source $S$, such that $W^{uu}(S)$ intersects $W^s(P)$.}  \label{fig:bicycle}
\end{center} 
\end{figure}

To overcome the lack of Lebesgue measure on $\Diff(M)$, several conjectures from the 90s used the concept of Kolmogorov typicality:
\begin{definition}\label{def: Kolmo typi} 
For $1\le r\le \infty$ and $d\ge 1$, a property $( P)$ is \emph{$d$-$C^r$-Kolmogorov typical} if for every compact manifold $\cP$ of dimension $d$, for a generic $C^r$-embedding  $p\in \cP \mapsto f_p\in \Diff(M)$,  the map $f_p$ satisfies the property $( P)$ at Lebesgue a.e.  parameter $p\in \cP$. 
\end{definition} 
In other words, a Kolmogorov typical property is true at Leb. a.e. every  parameter of $d$-dimensional generic family $(f_p)_{p\in \cP} $.  Using this notion, Pugh and Shub formulated the following conjecture:
\begin{conjecture}[{\cite[Conj. 3]{PS96}}] \label{PSconj}
For every $d\ge 1$, a $d$-Kolmogorov typical diffeomorphism displays finitely many sinks and other attractors.
\end{conjecture} 
This viewpoint was also promoted by Palis in one version of his global conjecture~\cite{Pa08}: ``In brief, for a typical dynamical system, almost all trajectories have only finitely many choices, of (transitive) attractors, where to accumulate upon in the future”.  The following is in opposition to  this vision and  disproved Pugh-Shub's conjecture in the finite regular case:
\begin{theorem}[\cite{Be16,Be17}] For every $1\le r< \infty$, for every $d\ge 0$, there exists a nonempty open subset $ \cU_d$ of $d$-dimensional $C^r$-families of $C^r$-local diffeomorphisms of the annulus, such that a generic
 family\footnote{This common way of writing means that  there exists a topologically generic subset $\cal G$ of $\cU_d$  such that for every $(f_p)_p\in  \cal G$, for every parameter $p$, the map $f_p$ displays  infinitely many sinks.} $(f_p)_p\in  \cU_d$ displays infinitely many sinks at every parameter $p$.

The analogous result is true for diffeomorphisms of  manifold $M$ of dimension $\ge 3$.
\end{theorem} 
\begin{proof}[Sketch of proof] By Newhouse's theorem,  an arbitrarily small  perturbation of a dissipative homoclinic tangency creates a sink.  Then
for every compact manifold $\cP$,   it is sufficient to find a locally  dense\footnote{A \emph{locally dense (resp. generic)  subset} is a subset which is dense (resp. generic) in a nonempty open subset.} subset of families $(f_p)_{p\in \cP}$ which displays a dissipative saddle point $\Omega$ 
such that for every $p_o\in \cP$, the saddle point $\Omega$  displays a homoclinic tangency which persists on a neighborhood of $p_o$. 
In order to do so, we assume that $\Omega$ displays a   stable manifold $W^s(\Omega ;f_p)$ which intersects the repulsion basin of a source $S_p$.  Hence it accumulates at the source.  
Then, we introduced a family of hyperbolic sets $(K_p)_p$ called \emph{parablender} which satisfies, roughly speaking, the following property. 
In the induced dynamics on the space of $C^d$-jet of families of points:
\[J^d_{p_o}  (z_p)_p\mapsto J^d_{p_o} (f_p(z_p))_p ,\]
the jets of the  continuation of points in  $(K_p)_p$ forms a blender. 
Hence there is a parablender whose local unstable lamination $W^u_{loc}(K)$ satisfies that for every $p_o$, there is $k\in K$,  
such that the distance from $S_{p}$ to  $W^u_{loc} (k; f_{p})$  is $o((p-p_o)^d)$.  Consequently, a $C^d$-perturbation of the family makes $S_p$ belong to $W^u_{loc} (k; f_p)$ for every $p$ in a uniform neighborhood of $p_o$. 
As $S_p$ is close to $W^s(\Omega ;f_p)$ for every $p$,  by connecting  $\Omega_p$ to $K_p$, the inclination lemma implies that $W^u_{loc}(\Omega ;f_p)$  and $W^s(\Omega ;f_p)$ are $C^d$-close to be tangent in a uniform neighborhood of $p_o$.   
From this we deduce the desired property.
\end{proof}
In~\cite{BCP22}, we strengthened part of the above theorem  by exhibiting  an open set $\cB\subset \Diff_{loc}^r(M^2)$ of surface local diffeomorphisms such that any $C^r$-generic one-parameter family  $(f_p)_{p\in  \R}$ in  $\cB$ is such that $f_p$ displays infinitely many sinks for every small $p$.   We called this notion ``germ typicality”.  Importantly, we describe the set $\cB$ as accumulated by the configuration depicted in \cref{fig:bicycle}, called {bicycle}. 

All the above results deals with the finitely regular family. A natural question asked many times by Yoccoz is whether these theorems are also valid for infinitely regular families. 
An ongoing collaboration with Biebler makes us confident in the following conjecture (which would answer to Yoccoz' question):
\begin{conjecture}
For every $1\le r\le \infty$, there is  a nonempty open set $\cB\subset \Diff_{loc}^r(M^2)$ in which the coexistence of infinitely many sinks is $d$-$C^r$-Kolmogorov typical for every $d\ge 0$.
The analogous result is true for diffeomorphisms of  manifold $M$ of dimension $\ge 3$.
\end{conjecture} 
We notice that Pugh-Shub's conjecture and Palis' vision remains strong in the category of  surface diffeomorphisms.  The formulation with parameter family of Palis' conjecture was motivated by the Benedicks-Carleson  breakthrough and subsequent works:
\begin{theorem}[\cite{BC2,BV01,BY00,MV93,WY08,Ta11,berhen}] \label{BC thm}
For every $b\in \R$ small enough,  there exists a neighborhood $ \cU$ of the family $(h_a)_a$, with $h_a:= (x,y)\mapsto (x^2+a-y, bx)$, such that for every $(\tilde h_a)_a\in  \cU$, there exists a parameters set of positive Lebesgue measure at which $\tilde h_a$  displays a stochastic attractor (SRB measure)  whose basin contains Leb. a.e. point which does not escape to infinity. \end{theorem} 
The latter theorem states the abundance of almost unique  ergodicity for  (strongly) dissipative Hénon-Like maps (which can be shown to lie within the Newhouse domain). 
In the opposite direction, we can wonder about the abundance of the infinitude of the attractors. This question goes back to the work of 
Tedeschini-Lalli and Yorke~\cite{TLY86}.
\begin{problem}\label{abundanceNewhouse}  Are there compact  manifolds $M$,  $\cP$ and an open set $ \cU$ of smooth families $f_\cP:=(f_p)_{p\in \cP} $ of maps $f_p\in C^\infty(M, M)$ such that the following subset of $\cP$ has positive Lebesgue measure:
\[ {\cQ}(f_\cP):=  \{ p\in \cP: f_p\text{ displays infinitely many attractors}\}\; ?\]
\end{problem} 
A positive answer to this problem would certainly bring about a revolution in many fields.  Indeed we would have a model for systems such that for many initial conditions, the  observables  converge in mean, but whenever we redo an experiment, this statistical mean would  change dramatically.
The same causes would not produce necessarily the same effects in the statistical  paradigm.  To study this problem, one can start by studying the Hausdorff dimension of Newhouse phenomenon. With  de Simoi~\cite{BS16}, we showed that among one-parameter family of surface diffeomorphisms, during the one-dimensional unfolding of a dissipative quadratic homoclinic tangency, the coexistence of infinitely many sinks appears on a subset of Hausdorff dimension $\ge 1/2$. See also~\cite{Wa90}. 
The following proposes to improve  this bound by enlarging the dimension of the spaces.
\begin{problem}  Are there   compact manifolds $M$,  $\cP$  and an open set $ \cU$ of smooth families $f_\cP:=(f_p)_{p\in \cP} $ of maps $f_p\in C^\infty(M, M)$ such that: 
\[ \dim_{HD} \cQ(f_\cP) >\dim \cP - \tfrac12\; ?\]
\end{problem}

\subsection{A symplectic counterpart of the problem of the finitude of attractors}
Let $(M,\omega)$ be a compact symplectic surface. 
Traditionally, given a  symplectomorphism in  $\Symp(M)$, the surface $M$  is split  into two regions: a chaotic part formed by sets called  stochastic seas  and a quasi-periodic part formed by domains called elliptic islands. 
Let us give an interpretation of this picture.  The stochastic sea is assumed to display an ergodic behavior while the  elliptic islands are usually imagined as disks, each surrounding an elliptic point and bounded by a KAM curve.
 Yet, each elliptic periodic point is typically accumulated by  infinitely many elliptic periodic points. Hence there are typically uncountably many cascades of  elliptic periodic points. Furthermore, each elliptic periodic  point is typically surrounded by a homoclinic  tangle, conjectured to lie inside   a stochastic sea.  In this generality, elliptic islands and seas are nested, and there are infinitely many of them. Yet, up to approximation, using Birkhoff normal form, we can assume that an open neighborhood of each elliptic point or KAM torus is integrable\footnote{In higher dimension, \emph{integrability} is taken here in a strong sense: there is a partition of the domain by invariant tori on which the dynamics acts as a rotation.}
:

\begin{definition}
An  \emph{integrable elliptic island} is a domain onto which the dynamics is integrable, and which is maximal with this property.  
A \emph{stochastic sea} is a subset of positive Lebesgue measure on which the dynamics is ergodic.
\end{definition} 
This enables us to wonder about a symplectic counterpart of  \cref{Density finite}:
\begin{problem}   
Is there a dense subset of $\Symp^\infty(M)$ formed by dynamics displaying finitely many integrable elliptic islands  $\bigsqcup I_j$ 
and finitely many stochastic seas  $\bigsqcup_i S_j$, such that the union  $\bigsqcup I_j\cup \bigsqcup_i S_j$ is Leb a.e. $M$?
\end{problem} 
\begin{remark} 
The above definition and problem make sense in any dimension.
\end{remark} 
 
\subsection{Strange attractors, stochastic seas and strong regularity}
The study of wild dynamics was motivated also by old conjectures based on numerical experiments.
First Lorenz  studied a certain model for the meteorology.
He discovered and depicted   the shape of  his eponymous   attractor arising from  a 3-dimensional flow and noticed that this dynamics is sensitive to the initial conditions from   numerical observations~\cite{Lor63} .  
Around the same time, Hénon and Heiles studied a Hamiltonian flow modeling galactic motion~\cite{HW64}. They performed numerical experiments by plotting  orbits on a 2-dimensional Poincaré  section and conjectured the existence of an open subset on which the dynamics is ergodic.  Similar observations were then done by Chirikov~\cite{Chi69} for iteration of the following  the so-called \emph{standard} family of symplectic maps of the torus:
 \[  S_K:  (\theta, p)\mapsto  (
\theta + p + K\cdot \sin(\theta) 
,
p+K\cdot \sin(\theta) )\; , \quad K\in \R.\]
Then Hénon~\cite{Hen76} showed  numerical evidences that the iteration of the following map  displays a ``strange attractor'': 
\[ (x,y)\mapsto (1- a\cdot x^2 +y, b\cdot x)\quad \text{with } a= 1.4 \text{ and } b= 0.3\;.\]
 Tucker~\cite{Tuc02} showed that the Lorenz flow satisfies a model  of Williams and that it is in particular partially hyperbolic: it expands uniformly  one direction, and this more than any other directions. Then the theory of partially hyperbolicity enables a tame understanding of the statistical behavior of this dynamics.  It is not the case of the examples of  Hénon–Heiles,  Chirikov and Hénon.

 The Hénon Conjecture is  open. The most promising work in this direction is certainly the Benedicks-Carleson \cref{BC thm}.  The initial proofs  of the latter theorem were given by defining a map whose fixed point  displays a local unstable  manifold whose iterates have a geometry which satisfies several analytic estimates. To formulate these conditions, they introduce the notion of critical point; a point  analytically defined along the induction as sent to the tip of a folded curve. 
 Then the existence of parameters in the Hénon family which satisfy such conditions is done using a parameter selection, similarly to Jakobson's Theorem   in one-dimensional dynamics~\cite{Ja81}. Yet a difference  is that the number of critical points increases exponentially fast with the induction, and when the parameter varies, some critical points may vanish. A tricky part of the proof is the description of  the parameter dependence of    this structure.
 
Yoccoz' \emph{strong regularity} program~\cite{Yolecture1,BY19} proposed to  reprove these theorems and go further, by using combinatorial and topological definitions, and deduce from them analytic estimates.  This paradigm has been successful  in
dynamics in  one-complex variable (puzzle piece).
 
\begin{figure}[h]
\begin{center} 
 \includegraphics[height=5.cm]{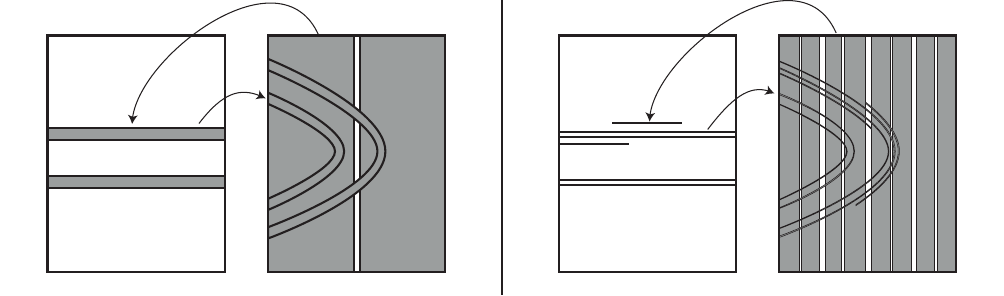} 
 \caption{
Each vertical box is sent by an iterate to a horizontal box and then  folded at another vertical box. First  and second steps of the induction. } \label{figSR}
\end{center} 
 \end{figure}

Yoccoz  wrote a  proof of Jakobson's theorem using this method~\cite{Y19}.  Then Palis-Yoccoz continued this program by studying the bifurcation of homoclinic tangencies of horseshoes with fractal dimension slightly greater than one. 
In~\cite{berhen},  we achieved another step of Yoccoz' program  by giving an alternative proof of Benedick-Carleson's theorem using a two-dimensional counterpart of Yoccoz' puzzle pieces.
The idea of this new proof is to define   inductively on $k\ge 0$, and following  only combinatorial and topological conditions,   a family $(B_\sg)_{\sg \in \sG_k}$ of long vertical  boxes 
which are sent by the dynamics to   horizontal boxes $B^\sg$   by a certain iterate of the dynamics $f^{n(\sg)}$. 
The definition of $k$-strongly regular map is formulated by asking each of these horizontal boxes $(B^\sg)_{\sg \in \sG_k}$ to be folded by  the dynamics $f$ \emph{inside} a certain vertical box.
See \cref{figSR}.
 Such a topological condition implies a tame control of the geometry of all the boxes of the next generations $(B_\sg)_{\sg \in \sG_{k+1}}$ and $(B^\sg)_{\sg \in \sG_{k+1}}$. At the limit  the decreasing intersections $\bigcap_k \bigcup_{\sg \in \sG_k}  B_\sg$
and  $\bigcap_k \bigcup_{\sg \in \sG_k}  B^\sg$  form a vertical and horizontal laminations which are   combinatorially defined local  stable and unstable  manifolds. The folding condition implies that each of these local unstable manifolds is sent by $f$ to a curve which displays a quadratic tangency with one of these local stable manifolds.  We are indeed dealing with a wild hyperbolic set. 
 
 This approach leads to the following toy model  of the parameter selection:  how to include a Cantor set $K^s+a $ into another $K^u$ for a set of translation parameters $a$ of positive Lebesgue measure. In~\cite{BM16}, we provided a criterion for this model, asking $K^u$ to have in particular a positive Lebesgue measure, and the Hausdorff dimension of $K^s$ to be smaller than 1 minus the exponential rate of decay of the size of the gaps of $K^u$. 
 When $b$ is small, this condition is satisfied, for the fractal dimension of $K^s$ is dominated by $1/|\log b|$.  Also, in \cref{BC thm}, the  fractal dimension of the attractors should be  (very) close to $1$.
   \begin{problem}[Yoccoz~\cite{Yolecture1,BY19}] Show the abundance of strongly regular attractors in the plane with fractal dimension higher and higher.\end{problem} 
    The  case of the standard map would correspond to a fractal dimension equal to $2$.

\section{Wild dynamics from renormalization}
Let $M$ be an $n$-dimensional manifold. 
There are many definitions of renormalization. In this section we will work mostly with the following:
\begin{definition}
A map $F: U\subset \R^n \to V\subset \R^n$ is a \emph{primitive renormalization} of $f\in \Diff(M)$ if there is a diffeomorphism $\phi: U\cup V\hookrightarrow M$ and an integer $N$ such that:
\[ f^N\circ \phi|_U = \phi \circ F\]
and the sets $(f^k\circ \phi(U))_{0\le k< n}$ are disjoint.
\end{definition}
Analogously, when we study dynamics given by a symplectomorphism $f$, we should require the map $\phi$ to be conformally symplectic so that $F$ is a symplectomorphism as well.

\subsection{Renormalization near homoclinic saddle point} 
\begin{figure}[h]
\begin{center} 
 \includegraphics[height=3.cm]{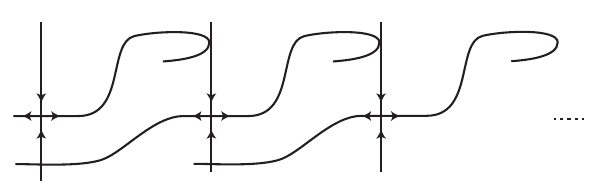} 
 \caption{Chain of heteroclinic tangencies between saddle points of a same horseshoe.} \label{ChaineTangence}
\end{center}
\end{figure}

For surface diffeomorphisms or symplectomorphisms, by the Newhouse and Duarte theorems~\cite{Ne79,Du08}, a small perturbation of a homoclinic tangency displays a wild hyperbolic horseshoe. From this, we can easily deduce a perturbation which displays an arbitrarily large number $N$ of distinct   quadratic heteroclinic tangencies between connected saddle points, as in \cref{ChaineTangence}. 
These $N$-heteroclinic tangencies can be unfolded independently along an $N$-dimensional parameter family. 
Then the  Gonchenko-Shilnikov-Turaev Theory~\cite{Tu03,GTS07,BFP23} allows  to deduce the following family of renormalizations:
 
\begin{theorem}\label{GST thm} 
If $\Omega$ is area contracting, then for every $\epsilon>0$, there exists a domain $I_\epsilon$ near $0\in \R^N$ and a primitive renormalization  of $(f_p)_{p\in I_\epsilon}$ which is of the form:
\[F_p: (X,Y)\mapsto (R_p(X) -Y, b_p \cdot X)+E_p(X,Y)\; ,\]
where the $C^r$-norm of $(E_p)_p$ is  $<\epsilon$ for every $r<1/\epsilon$ and $p\in I_\epsilon \mapsto  R_p \in \R_{N-1} [X]$ is a diffeomorphism onto the ball of radius $1/\epsilon$, and $(b_p)_p$ is small.   

Also if $f_p$ is a symplectic map and so $\Omega_p$ has determinant $1$, the same is true with $b_p=1$ for every $p$. 
\end{theorem} 

This result is very useful; for instance, it yields a short alternative proof of a theorem of Kaloshin.
\begin{theorem}[\cite{Ka00}] \label{KaloThm}For every surface $M$, 
for every sequence $(n_k)_{k\ge 0} $ of integers, there exists a topologically generic subset $\cG$ in the  Newhouse domain of $\Diff(M)$ such that any $f\in \cG$ displays at least  $n_k$ hyperbolic points of period $k$ for infinitely many $k$:
\[ \limsup \frac1{n_k} Card\,  \{\text{ hyperbolic fixed points }  x\in M  \text{ of }f^k\}\ge 1\; .\]
In particular, for every surface $M$, there is no topologically generic subset of $\Diff(M)$ which is formed by diffeomorphisms 
whose number of periodic points grows at most exponentially fast.
\end{theorem}
\begin{proof}
We start with any diffeomorphism $f$ in the Newhouse domain. After a perturbation as above, for a large  integer $N$, we obtain a family of perturbations which displays a family of renormalizations $(F_p)_p$ as in \cref{GST thm}, with $(b_p)_p$ small (up to 
reversing the dynamics). We obtain a parameter $p\in I_\epsilon$ at which the map $F_p$ displays a semi-parabolic fixed point near $0$ with 
a strong stable direction and a central local manifold  
on which the next $N-1$ derivatives of the restriction vanish at the fixed point.  
By a $C^N$-small perturbation (changing only the $C^{N+1}$-jet of the dynamics at this point),  we can assume it moreover   topologically unstable along its central manifold. Then its central manifold is smooth nearby this point. 
Hence after an arbitrarily $C^N$-small and smooth perturbation of $F_p$ and so $f_p$,  we obtain a normally contracted curve on the surface along which the dynamics is the identity. Therefore one can perturb the dynamics to produce more than $b_k$ saddle periodic points for some $k$ large. Such points persist in an open set of perturbations.
\end{proof} 
 \begin{remark} Kaloshin  proved the local genericity in   $\Diff^r(M)$ for  $2\le r<\infty$, but not for $r=\infty$ as above.  
     \end{remark} 
 
 In~\cite{BerPer}, we showed the counterpart of \cref{KaloThm}   for every manifold of dimension $\ge 3$, using a different methods (involving KAM). Moreover it addressed a problem posed by Arnold from the 90s asking to 
  show the Kolmogorov-typicality of diffeomorphisms whose number of periodic points grows at most exponentially fast. More precisely we showed the opposite answer in the finite regular case:
\begin{theorem}[\cite{BerPer}]
For every sequence $(n_k)_{k\ge 0} $ of integers, for every $1\le r< \infty$, for every $d\ge 0$, there exists a nonempty open subset $ \cU_d$ of $d$-dimensional $C^r$-families of $C^r$-local diffeomorphisms of the annulus, such that a generic family $(f_p)_p\in  \cU_d$ satisfies that for every $p$,  the map $f_p$ displays more than $n_k$ hyperbolic points of period $k$ for infinitely many $k$.

The analogous result is true for diffeomorphisms of manifolds $M$ of dimension $\ge 3$.
\end{theorem} 
\begin{proof}[Sketch of proof]
First observe that if a map of a circle displays a unique fixed point which is parabolic, then it can be perturbed to be conjugated to an irrational rotation, and (using KAM) this rotation can be perturbed in turn to a rational rotation. This enables to create as many periodic saddle points as desired for some large periods. We prove the counterpart of this result for parameter families and we then use a construction involving a parablender to show the existence of a locally dense set of families displaying at every parameter a  
normally hyperbolic circle whose inner dynamics remain parabolic under perturbations along the family. 
\end{proof}
\subsection{Universality} 
The notion of universality was introduced in $C^1$-dynamics by Bonatti and D\'iaz~\cite{BD03}, and then adapted by Turaev~\cite{Tu03} to smooth dynamics.  A diffeomorphism $f$ is $C^r$-\emph{universal} if the set of its primitive renormalizations is dense in the space of orientation preserving $C^r$-embeddings of the $n$-ball $\B$ into $\R^n$. Similarly, a symplectomorphism is $C^r$-\emph{universal} if the set of its primitive renormalizations is dense in  the space of $C^r$-symplectic embeddings of the $n$-ball $\B$ into $\R^n$. We notice that universality is a $G_\delta$-property.  
Bonatti and D\'iaz  showed:
 
\begin{theorem}[\cite{BD03}] \label{BD03} For every compact manifold $M$ of dimension $\ge 3$, $C^1$-universal dynamics form a locally topologically generic subset of $ \Diff^1(M)$.
\end{theorem}
Our next theorem generalizes one in dimension 2 by Gelfreich and Turaev~\cite{GeT10}:
 \begin{theorem}[\cite{BT25b}]\label{BT25b} 
 Let $(M, \omega)$ be a symplectic manifold (of any dimension) and let $\Symp^\infty_{e }(M)$ be the space of $C^r$-symplectomorphisms of $M$ which display a totally elliptic point\footnote{A \emph{totally elliptic point} is a  periodic point whose eigenvalues are all of modulus 1 but none of them is 1.}. Then 
there is a topologically generic subset $\Symp^\infty_{e }(M)$ formed by universal symplectomorphisms. 
 \end{theorem}
 For non-conservative surface dynamics, Turaev showed the following:
\begin{theorem}[\cite{Tu15}] \label{Tu03}For every  surface $M$ and $\infty\ge r\ge 2$, $C^r$-universal dynamics form a locally generic subset of $ \Diff^r(M)$.
\end{theorem} 
The proofs of these theorems are done by showing first the existence of dense subsets of dynamics  displaying  a periodic spot:
\begin{definition} A \emph{periodic spot}   is a nonempty open subset formed by periodic points. \end{definition}
\begin{proof}[Sketch of proof of \cref{BD03}] Using a blender, Bonatti and D\'iaz  construct a wild hyperbolic set such that a dense subset of its $C^1$-perturbations is formed by dynamics  displaying periodic points whose eigenvalues are all equal to $1$. Such periodic points can be $C ^1$-perturbed to lie into a periodic spot. Then we achieve the proof by invoking the Ruelle-Takens theorem~\cite{RT71}
which implies that any direct embedding $\B\hookrightarrow M$ can be $C^1$-approximated by a primitive renormalization of a perturbation of the identity. \end{proof}
\begin{proof}[Sketch of proof of \cref{BT25b}] 
We first perturb the elliptic point to make it integrable and non-degenerate. Then it is accumulated by ``integrable” parabolic periodic points. 
Consequently there are  primitive renormalizations of  neighborhoods of them in which the symplectomorphism is a linear parabolic map. This can be easily perturbed to a rational rotation. This gives the density of  symplectomorphisms displaying periodic spots. 

We can perturb these spots to the time-$\epsilon$ of a product of pendulums flows. Such displays a flat homoclinic tangency. Then in dimension 2, we can apply \cref{GST thm} which produces renormalizations approximating any generalized Hénon maps. Then it suffices to recall that any symplectic map can be approximated by a primitive renormalization of a generalized Hénon map~\cite{Tu03, BT25a}. 
 \end{proof}
\begin{proof}[Sketch of proof of \cref{Tu03}] 
Using a wild hyperbolic horseshoe, with both area contracting and expanding periodic points, Turaev showed the existence of a locally dense set of perturbations exhibiting periodic spots. Then
he strengthened Ruelle-Taken's Theorem by proving that any direct $C^\infty $-embedding $G$ of $\B^n$ into $\R^n$  can be {\em approximated} by a primitive renormalization of a $C^\infty$-perturbation $g$ of the identity. 
 \end{proof} 
With Helfter and Gourmelon, we proved the following improvement of the second part of the latter  argument (the approximation becomes an equality):
\begin{theorem}[\cite{BGH24}]\label{BGH24}
For any $1\le r \le \infty$ and any orientation preserving $G\in \Diff^r(\B^n)$, in any neighborhood $N\subset \Diff^r_c (\B^n)$ of the identity, there exists $g\in N$ such that a primitive renormalization
of $g$ is {\em equal} to $G$. 
\end{theorem}
\begin{question}
Given any $G\in \Symp^\infty(\B^n)$ and any neighborhood $N\subset \Symp^\infty (\B^n)$ of the identity, is there $g\in N$ such that a primitive renormalization
of $g$ is equal to $G$? 
\end{question}
We remark that the following problem remains open:
\begin{problem} Show that for any manifold $M$ of dimension $n\ge3$ and $2\le r\le \infty$, there exists a locally generic set of $\Diff^r(M)$ formed by $C^r$-universal maps. 
\end{problem} 
To solve this problem, in view of~\cite{Tu03,BGH24}, it suffices to answer:
\begin{problem} Show that for any manifold $M$ of dimension $n\ge3$, there exists a nonempty locally dense subset of $\Diff^\infty(M)$ formed by diffeomorphisms displaying a periodic spot.
\end{problem}
In~\cite{BFP23}, we showed the existence of  steady Euler flows of $\R^3$ exhibiting a wild hyperbolic sets. Then given a  generic such a flow, the set  its Poincaré return maps  form a dense subset of $\Symp^\infty(\D\hookrightarrow \R^2)$. 
Also  Turaev  communicated to me  the following: 
\begin{question}
Does universality exist in the category  of analytic surface diffeomorphisms? How about the category of conservative polynomial automorphisms?\end{question}

\subsection{Density of positive metric entropy} 
Herman's positive entropy conjecture~\cite{He98} asserted the existence of a $\Symp^\infty$-perturbation $f$ of the identity on the disk exhibiting positive metric entropy $h_\leb(f)$:
\[ h_\leb(f):= \lim_{n\to \infty} \int \frac1n \log \|Df^n\| \, d\leb\; .\]
 We established the following stronger result:
\begin{theorem}[\cite{BT19}] 
\label{BT19}
If $f$ is a $C^\infty$-symplectomorphism displaying an elliptic point, then it can be $C^\infty$-perturbed to one with positive metric entropy.  
\end{theorem}
 \begin{figure}[h] \begin{center} 
 \includegraphics[height=4.cm]{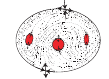} 
 \caption{A stochastic sea given by the surgery of a linear hyperbolic map.} \label{stocha_sea}
\end{center} 
\end{figure}
\begin{proof}[Sketch of proof] By    \cref{BT25b}, the renormalizations of a perturbation of an elliptic point approximate any dynamics. We chose   one  which displays a stochastic sea, as depicted  \cref{stocha_sea}, given by a surgery (blowup of 4 fixed points) of a linear hyperbolic map of the torus. 
Yet the approximation may not have positive metric entropy, but it is still the case if the heteroclinic connection persist (coincidence of half stable and unstable manifolds). To obtain such,  
we show that the renormalization can be equal to this  approximation composed with any  vertical shear $(x,y)\mapsto (x,y+g(x))$. The function $g$ is then found by studying the Melnikov operator in neat coordinates.
\end{proof} 
The first step of the latter argument is constructed in~\cite{Ber23} in the analytic category. 
The following remain open:
\begin{problem} 
 Show the existence of an analytic perturbation of the identity on the disk with positive metric entropy.
\end{problem}
\begin{problem} Show that a perturbation of an ellipse  defines a billiard with positive metric entropy. \end{problem}
\begin{problem} 
 Show that an integrable surface geodesic flow can be perturbed to display a positive metric entropy. 
\end{problem} 
The proof \cref{BT19} might be useful to solve these problems. However the heteroclinic  connections  appearing in this argument  are very fragile; also they cannot appear among entire symplectomorphisms~\cite{Us80}, such as the conservative Hénon map or the Chirikov standard map. 

For any symplectic surface $(M,\Omega)$, 
either {\em (a)}  zero metric entropy is a locally generic property in $\Sym^\infty(M)$, or 
{\em (b)}  positive metric entropy is an open-dense property in  $\Sym^\infty(M)$.
\begin{question}[Compare with {\cite[Q. 4.3]{He98}}] Is (a) or (b) true ?\end{question}
By \cref{BT19}, to show the density of positive metric entropy,  it   suffices to show that if a surface symplectomorphism is weak-$*$ stable (none of its perturbations displays an elliptic point), then it is uniformly hyperbolic;  a  symplectic counterpart of a conjecture by Ma\~ne.

\subsection{Generic local property and Kolmogorov typicality} 
Let $\cF(M)$ be a Fréchet manifold formed by dynamics such as  $\Diff^r(M)$, $\Symp^r(M)$, for $1\le r\le \infty$, or open subsets of these spaces. 
We would like to study what are the properties which are typical in the following sense.  
\begin{definition} A property $( P)$  is \emph{Kolmogorov typical} in $\cF(M)$ if it is $d$-Kolmogorov typical  for every $d\ge 0$ (see def. \ref{def: Kolmo typi}). 
The property $(P)$ is \emph{locally Kolmogorov typical} if there exists  a nonempty open subset $\cV\subset \cF(M)$  in which $(P)$ is Kolmogorov typical.  
 \end{definition}
Note that any property which is  Kolmogorov typical is  topologically generic   (take $d=0$).  
In order to study the converse, we introduced in~\cite{Ber19b} the following notion:

\begin{definition}\label{localizable equiv prop}
A property $(P)$ is \emph{localizable} if, for any nonempty open set $U\subset  M$,  there exists a map close to the identity, compactly supported $g\in \Symp^\infty (U)$ 
such that $(P)$ holds for any $f\in \Symp^\infty(M)$ satisfying $f^N|_U=g$ for some $N\geq 1$.  
It is   \emph{openly localizable} if $g$ can be chosen so that the property $(P)$ is satisfied by any map in a neighborhood in $\Symp^\infty(M)$  of the set of maps $f$ coinciding with $g$.  
\end{definition} 

For instance, exhibiting a wild hyperbolic set is an openly localizable property in any $\Diff(M)$, for manifolds $M$ of dimension $\ge 2$. 
\begin{definition} 
A property $(P)$ is \emph{generically localizable} if it is implied by the conjunction of a countable family of openly localizable properties $\{(P_i): {i\ge 0}\}$:
\[\bigwedge_i (P_i)\Rightarrow (P)\; .\]
\end{definition} 

For instance, by \cref{BT25b,BGH24}, universality is a generically localizable property  
in $\Symp^\infty(M)$ and in $\Diff^r(M)$ for every $\infty\ge r\ge 2$.  
 Note also that the super-exponential growth of the number of periodic points is a generically localizable property (see~\cite[Prop. 2.2]{Ber19b}). 

We can now state:
\begin{conjecture}[\cite{Ber19b}]\label{Conjecture}
There are many space of dynamical systems for which any generically  localizable property is locally Kolmogorov typical.
\end{conjecture}
This conjecture has been proved 
for the space of symplectomorphisms:
\begin{theorem}[\cite{BT25b}]  \label{BT25c}
For every symplectic manifold $(M,\Omega)$, any generic localizable property is  Kolmogorov typical in the open set $\Sym^\infty_{e\ell} (M)$ of symplectomorphisms possessing a non-degenerate totally elliptic point.
In particular  universality is   Kolmogorov typical in $\Sym^\infty_{e\ell} (M)$.  
 The same occurs with the super-exponential growth of  the number of periodic points.
  \end{theorem} 
 The Kolmogorov typicality of the super-exponential growth of  the number of periodic points was proved for surface symplectomorphisms by  Asaoka~\cite{Asa17}.
 \section{Emergence}
When a dynamical system on a compact manifold $M$ preserves the Lebesgue measure, Birkhoff’s ergodic theorem implies that the following sequence converges for Lebesgue almost every $x \in M$: 
\[ \se_n(x) := \frac1n \sum_{0\le k<n }\delta_{f^k(x)}\;  \] 
to a probability measure $\se(x)$ called the \emph{empirical measure} of $x$.  The empirical measure describes  
the statistical behavior of the orbit of $x$.  The distribution of these statistics $\{\se(x)\, :\, x\in M\}$  is given by the pushforward $\se_\star \Leb$ of the Lebesgue measure onto the space of probability measures $\cM(M)$.  The measure  $\se_\star \Leb$ is a measure on the space of measures. It is called the \emph{ergodic decomposition}.  

When the Lebesgue measure is ergodic, then $\se(x)=\Leb$ a.e. Hence the ergodic decomposition is a Dirac measure at $\Leb$. 
If the dynamics is almost finitely ergodic, then there exists a finite set of  probability measures $\{\mu_1, \dots, \mu_k\}$  to which  $\se(x)$ belongs for almost every $x\in M$.  Then the ergodic decomposition is a finite sum of atoms $\sum_i \lambda_i \delta_{\mu_i}$.  Note that the identity is not finitely ergodic: its ergodic decomposition is $\int_M \delta_x \, d\Leb $. However, it is ``finite-dimensional”.  
 
In contrast to  Boltzmann’s ergodic hypothesis, we may ask when the ergodic decomposition is far from finite-dimensional.  
To quantify this, we introduced the notion of emergence~\cite{Be17}. Let us recall its definition. We endow the space of probability measures with the Kantorovich-Wasserstein distance:
\[ d: (\mu_1, \mu_2)\in \cM(M)\mapsto \sup_{\phi \in \mathrm{Lip}^1(M)} \int \phi\,  d(\mu_1-\mu_2)\; . \]
We recall that this distance endows  
$\cM(M)$ with the weak-$\star$-topology. 
\begin{definition} The \emph{emergence} $\cE(\epsilon)$ at scale $\epsilon>0$ is the minimal number $N$ of probability measures $(\mu_i)_{1\le i \le N}$  such that:
\[ \int_{x\in M} \min_{1\le i \le N} d(\se(x), \mu_i) \, d\Leb<\epsilon\; .\]
\end{definition}  
We say that the emergence is \emph{high} if  $\cE(\epsilon)$ is not dominated by a power of $\frac1\epsilon$:
\begin{equation} \tag{High Emergence} \limsup_{\epsilon\to 0} \frac{ \log \cE(\epsilon)  }
{|\log \epsilon|} =\infty\end{equation} 
Loosely speaking, the emergence is high when  the ergodic decomposition is infinite-dimensional. This is very far from the paradigm of the ergodic hypothesis.  

There exist non-conservative dynamical systems for which the sequence of empirical measures $(\se_n(x))_n$ fails to converge for Lebesgue almost every $x \in M$.
We do not know if this situation is negligible among smooth dynamics. In particular, the following is open:
\begin{problem}[Thom~\cite{Su74}, Ruelle, Takens, Pugh, Shub 70's-90s'] 
Does a generic (or typical) diffeomorphism in $\Diff(M)$ have its sequence of empirical measures $(\se_n(x))_n$ converging for Lebesgue almost every $x \in M$?
\end{problem}
Nevertheless we can define the emergence of a non-conservative map $f$ as: 
\begin{definition} The \emph{emergence} $\cE(\epsilon)$ at scale $\epsilon>0$ is the minimal number $N$ of probability measures $(\mu_i)_{1\le i \le N}$  such that:
    \[ \limsup_{n\to \infty}  \int_{x\in M} \min_{1\le i \le N} d(\se_n(x), \mu_i)\,  d\Leb<\epsilon\; .\]
\end{definition}  
It is easy to see that the emergence is at most the covering number $\cN_B(\epsilon)$ by $\epsilon$-balls of the space of probability measures.  
The latter displays the following asymptotic behavior when $\epsilon\to 0$:
\[ 
\text{ord }\cN_B:= \frac{\log \log \cN_B(\epsilon)}{|\log  \epsilon |}\sim \dim M\; . 
\]
Hence we say that the map $f$ has \emph{emergence of maximal order} if:
\[ \limsup_{\epsilon\to 0} \frac{\log \log \cE(\epsilon)}{|\log  \epsilon |} = \dim M\; .\]
This means that for every $\delta>0$, at infinitely many scales $\epsilon>0$,  one needs more than $\exp (\epsilon^{-(\dim M-\delta)})$ % added parentheses, fixed exponent
to describe the statistical behavior of the dynamics to precision $\epsilon$. The growth is super-polynomial, and thus cannot be feasibly computed by any algorithm. 
Note that  when $\dim M\ge 2$, this number is super-exponential. In~\cite{Be17}, I proposed the following conjecture:
\begin{conjecture}[\cite{Be17}] \label{Thm Typi emergence}
In many spaces of dynamical systems,   high emergence is  locally typical, in many senses.  
\end{conjecture}
This conjecture motivates several works, some of which solved conjectures of independent interest; we will recall   some of them in the next subsections. Similarly to \cref{abundanceNewhouse}, 
we can wonder  about the abundance of high emergence:
%a strong sense of typicality is wondered here:
\begin{problem}\label{abundanceemergence}  Are there compact  manifolds $M$,  $\cP$ and an open set $ \cU$ of smooth families $ (f_p)_{p\in \cP} $ of selfmaps $f_p$ of $M$  such that the following subset of $\cP$ has positive Lebesgue measure:
\[  \{ p\in \cP: f_p\text{ displays high emergence}\}\; ?\]
\end{problem}

\subsection{Genericity and typicality of emergence of maximal order} 
In~\cite{BB21}, with Bochi, we proved that a generic symplectomorphism displaying an elliptic point has emergence of maximal order. 
In~\cite{BT25b}, with Turaev, we generalized this local construction to any dimension,  applied 
\cref{BT25c} and  obtained the  following:
\begin{theorem}[\cite{BB21,BT25b}]\label{Thm Typi emergence sympl}
Kolmogorov typically in $\Symp_{e\ell}^\infty (M)$, a symplectomorphism displays an emergence of maximal order $\dim M$.
\end{theorem}
\begin{proof}By \cref{BT25c},  it suffices to show that exhibiting emergence of maximal order is a generically localizable property.   More precisely, it suffices to show that a perturbation of the identity has robustly high emergence at small scales. To this end, we embedded a twist map of $\T^n\times (0,1)^n $ in a highly distorted way, so that the 
tori which support the empirical measures are in means very apart of each other. By flattening  the twist  where the distortion blows up, this embedded dynamics can be extended smoothly. 
In this way, we obtain a perturbation of the identity with high emergence at arbitrarily small scales.
Finally,   we use  the KAM theorem to show that, for every small scale, there exists an open set of perturbations on which the emergence remains high.  
\end{proof}
Note that \cref{Thm Typi emergence sympl} confirms \cref{Thm Typi emergence}.  In~\cite{BB21}, we also proved  that displaying high emergence is a locally generic property among (nonconservative) surface diffeomorphisms. A consequence of  \cref{Conjecture} would then be that it is also locally Kolmogorov typical among diffeomorphisms of manifolds of sufficiently large dimension. 
\subsection{Maximal oscillation} 
In~\cite{HK95}, Hoefboer and Keller introduced the notion of \emph{maximal oscillation}: a self-map $f$ of a manifold $M$ is said to have maximal oscillation if, for Lebesgue  almost every $x$, the sequence of empirical measures $(\se_n(x))_{n\ge 0}$ accumulates onto the space of {\em all} invariant probability measures.  
Then they proved the existence of  unimodal maps $x\mapsto x^2+a$ which display maximal oscillation.  
By showing that such unimodal maps display a hyperbolic Cantor set of dimension arbitrarily close to $1$, and then recalling that its set of invariant measures has order close to $1$~\cite{BB21}, we obtain easily: % fixed grammar, capitalization
\begin{theorem} 
There exist unimodal maps $x\mapsto x^2+a$ with emergence order equal to $1$.\end{theorem} 
In~\cite{Ta22}, Talebi showed that rational functions on the sphere with maximal oscillation exist  (and are generic in the bifurcation locus). From this we deduce likewise:
\begin{theorem} 
There exist rational maps of the sphere with emergence order equal to $2$.\end{theorem}  

\subsection{Wandering stable component} 
Let $f$ be a self-map of a manifold.  A point $x$ is \emph{asymptotically stable} if it has a neighborhood $U$ formed by asymptotic points: for every $y\in U$, we have $d(f^n(x), f^n(y))\to 0$ as $n\to \infty$.  
A \emph{stable domain} is a connected open subset formed by asymptotically stable points.  A \emph{stable component} is a stable domain which is maximal.  
In other words, a stable component is a component of the set of asymptotically stable points. A stable component is \emph{wandering} if it does not intersect its iterates.  
 Together with Biebler, we proved the following result:
\begin{theorem}[\cite{BB23}]
There exists a locally dense set of real Hénon maps of degree 6, each exhibiting a wandering stable component $C$.
Moreover, the empirical measures of the points in $C$ accumulate on a set of invariant measures of positive order; in particular the emergence order is positive.   Furthermore, $C$ is the real trace of a wandering Fatou component.  
\end{theorem} 
The last part of this theorem answers a question posed by Friedland-Milnor~\cite{FL89}, among others, on the existence of wandering Fatou component for Hénon maps.  We recall that the Fatou set is the set of points which are Lyapunov stable (two points close have forward orbits that remain close).   See~\cite{Du22,BB23} which include surveys on the topic of wandering Fatou components in holomorphic dynamics. 

To prove this theorem, we show that if, within a parameter family, five quadratic homoclinic tangencies can be unfolded independently, and  are linked as in \cref{ChaineTangence}, then the estimate from~\cite{Ber18} can be used to implement an infinite chain of renormalizations, each of which is contracting. This implies the existence of an arbitrarily small parameter values of the unfolding that yield  wandering stable components.  The same argument  also proves the following:  

\begin{theorem}[\cite{BB23}]
For every $r\in \{2, \dots, \infty\} \cup \{\omega\}$, in the Newhouse domain, there exists a $C^r$-dense set of maps displaying a wandering stable component $C$.  
Moreover, the empirical measures of the points in $C$ accumulate onto a set of invariant measures of positive order.
In particular the emergence order of the dynamics is positive. 
\end{theorem} 

This theorem resolves the so-called Takens' last problem for $r=\infty$ and $\omega$.  For $r<\infty$, see the works of Colli, Vargas, Kiriki, Soma, and Nakano~\cite{CV01, KS17, KNS22}.  
\subsection{Pseudo-rotation}
If the dynamical system  $f$ is conservative, then its ergodic decomposition $\hat \se:= \se_\star \Leb$ is well defined. 
Recall that $\hat \se$  is a probability measure (on the space of probability measures $ \cM(M)$).  
Hence, by analogy with the notion of local dimension, we can define its local order by
\[ \text{ord }\hat \se (\nu):= \limsup_{\epsilon\to 0} \frac{|\log  \hat \se (B(\nu, \epsilon))|}{\log |\log \epsilon|}\quad  \text{ at any  $\nu\in \cM(M)$}.\]  

\begin{definition} The \emph{local emergence order} of $f$ is equal to 
$
\text{ord }\cE_{loc} = \int  \text{ord } \hat \se (\nu) \, d\hat \se (\nu)\;.  
$
\end{definition}

Helfter proved that this invariant is finer than the notion of emergence:
\begin{theorem}[\cite{He25}] The local emergence order is at most the emergence order:
\[ \text{ord }\cE_{loc}  \le \text{ord }\cE\;. \] 
\end{theorem}

In~\cite{Ber22}, using the AbC method of Anosov-Katok, we obtained the following.
\begin{theorem} 
There exists an analytic symplectomorphism of the cylinder $\bS^1 \times [0,1]$ whose local emergence order is maximal, that is $ \text{ord }\cE_{loc}=2$. 
\end{theorem}

Importantly, the AbC method was here implemented among \emph{analytic} symplectomorphisms. This development enabled the proof of the following:
\begin{theorem}[\cite{Ber22,Ber24}]
There exist transitive, analytic symplectomorphisms on the sphere, the disk, and the cylinder.  
\end{theorem}

This theorem resolves conjectures of Birkhoff~\cite{Bi27,Bi41} and  problems of Herman, Fayad-Katok  and Fayad-Krikorian~\cite{He98,FK04,FK18}.  
\begin{center}
\includegraphics[height=6.cm]{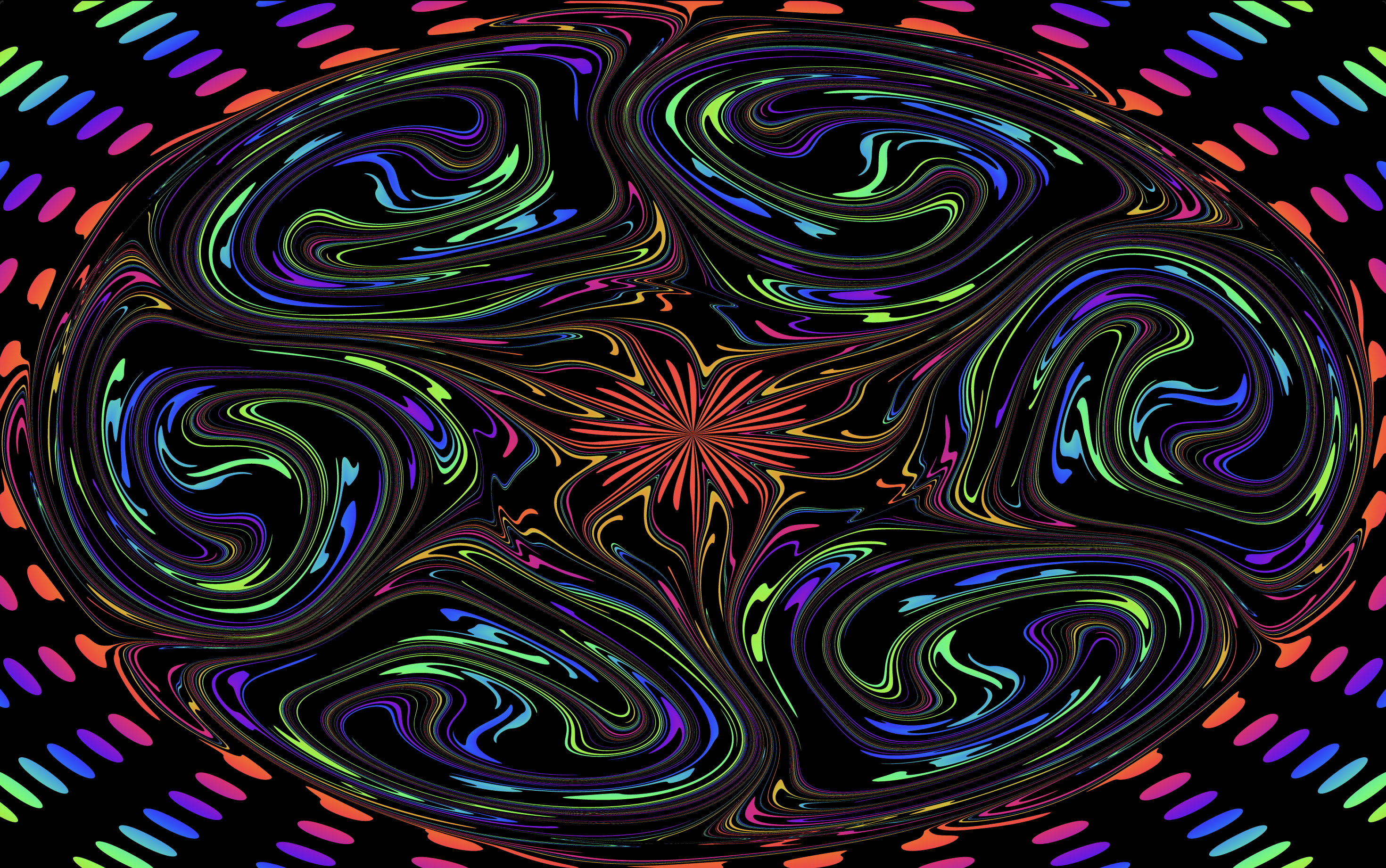} \\
{First steps of the AbC method on the disk.}
\end{center}
\begin{proof}[Sketch of proof] 
The AbC method constructs   symplectomorphisms of the form $f = \lim f_n$ with:
\[
f_n := h_1 \circ \cdots \circ h_n  \circ \text{Rot}_{\alpha_n} \circ h_n^{-1} \circ \cdots \circ h_1^{-1}\;, 
\] 
where $h_n$ commutes with $\text{Rot}_{\alpha_{n-1}}$ and $\alpha_n$ is a rational number sufficiently close to $\alpha_{n-1}$ to ensure that  $f_n$ is close to $f_{n-1}$.  
This implies that the limit $f = \lim f_n$ exists. To obtain $f$ transitive, the idea is to take $h_n$ that sends 
a circular orbit $(\text{Rot}_{\theta}(x))_{\theta\in \bS^1}$ to an orbit of $(h_{n-1} \circ \text{Rot}_{\theta}\circ h_{n-1}^{-1})_{\theta\in \bS^1}$
  which is $2^{-n}/L_n$ dense, with $L_n$ the Lipschitz constant of $h_{n-1} \circ \cdots \circ h_1$. 
  By taking the denominator of $\alpha_n$ sufficiently large, this implies that $f_n$ has a $2^{-n}$-dense orbit.
  Constructing $f$ sufficiently close to $f_n$ ensures that this property persists up to a small error, thereby implying the transitivity of $f$.
  The main difficulty in adapting this method to the analytic setting is to construct analytic maps that commute with a rational rotation of large denominator, yet significantly distort a circular orbit, and whose “radius of convergence” remains sufficiently large so that the composition  $h_n \circ \cdots \circ h_1$ has a ``convergence radius'' that does not converge to zero.  

To show this, we deform  the space and then we invoke the rigidity of analytic symplectic surfaces. For instance, in the case of the cylinder
$\A:=  \bS^1 \times [-1,1]$, the maps $h_n$ are taken as compositions of commutators of the form $[g_v, g_h]$, where $g_v(\theta,y) = (\theta, y + \psi(\theta))$ and $g_h(\theta,y) = (\theta + \phi(y), y)$, with $\psi$ an $\epsilon_n$-small trigonometric polynomial and $\phi$ a polynomial that is small on a set of the form $([-1/\epsilon_n, -1+2\epsilon_n] \sqcup [1-2\epsilon_n, 1/\epsilon_n]) \times [-i/\epsilon_n, i/\epsilon_n]$, with $\eta_n\ll \epsilon_n$ small.  
Such compositions are entire and are close to the identity near $\partial \A$.  
Hence, for sufficiently small $(\eta_n)_n$,  the sequence of compositions $ h_1 \circ \cdots \circ h_n$ restricted to the boundary $\partial \A$ converges to an analytic embedding of $\partial \A$ into $\bS^1 \times \R$.  
We thus obtain a transitive analytic symplectomorphism of a cylinder bounded by two analytic curves. Then it is easy to conjugate  this symplectomorphism  analytically to one which leaves invariant  $\A$.

This carries the case of the cylinder. For the disk and the sphere, using a bump function, we perturb these maps $h_n$ to maps $\tilde h_n$ which are symplectic, have a large complex extension that are {\em close} to being holomorphic, and whose restriction to $\C/\Z \times \{ \mathrm{Im}(w): |\mathrm{Re}(w)| \ge 1\}$ is the identity.  
Then the map
\[
\tilde f_n := \tilde h_1  \circ \cdots \circ \tilde h_n  \circ \text{Rot}_{\alpha_n} \circ \tilde h_n ^{-1}\circ \cdots \circ \tilde h_1^{-1}
\]
turns out to be analytic and symplectic, for the analytic structure given by the pushforward  by $\tilde h_1 \circ \cdots \circ \tilde h_n$ of the canonical one.  
Using a holomorphic extension of the symplectic polar coordinates $\D \to\A$ or of the axial projection of the sphere $\bS^2 \to\A$, we obtain a sequence of symplectomorphisms on the disk or the cylinder for an ``exotic'' analytic structure on the disk. Finally, we invoke the Newlander-Nirenberg theorem, which implies that these analytic structures converge.  
\end{proof}

Let us conclude by recalling that the AbC method produces many dynamics with interesting properties; see~\cite{FK04} for a survey. In~\cite{Ber24}, an AbC principle is introduced to import most of these properties from the differentiable setting to the analytic setting. For instance:
\begin{theorem}[\cite{De25}]
There exist analytic symplectomorphisms of the disk and the sphere whose local emergence order is maximal (=2). 
\end{theorem}

\section*{Acknowledgments.}
 I  thanks my masters, collaborators and students for sharing mathematics. In particular, J.-C. Yoccoz, M. Lyubich, J. Palis,   M.C. Arnaud, S. Crovisier, R. Dujardin, E. Pujals,  and D. Turaev. 
 
% SIAM recommends using BibTeX
% if using BibTeX
\bibliographystyle{siamplain}
%\bibliography{example_references}

\end{document}